\documentclass{amsproc} %orignal

\usepackage{amsmath,amssymb}
\newtheorem{theorem}{Theorem}[section]

\theoremstyle{definition}
\newtheorem{defi}[theorem]{Definition}

\theoremstyle{remark}
\newtheorem{remark}[theorem]{Remark}

\numberwithin{equation}{section}

\def\frak{\mathfrak}

\def\om{\omega}

\def\om{\omega}

\newfont{\df}{eufm10}
\def\vep{\varepsilon}

\def\ot{\otimes}

\def\ot{\otimes}

\def\om{\omega}
\def\ot{\otimes}

\begin{document}

\title[Fermionic
realization of $U_{r,s}(\widehat{\frak {sl}_n})$]
    {Fermionic
realization of two-parameter quantum affine algebra
$U_{r,s}(\widehat{\frak {sl}_n})$ }

%%% ----------------------------------------------------------------------
\author[Jing]{Naihuan Jing}
\address{Department of Mathematics,
   North Carolina State University,
   Raleigh, NC 27695-8205, USA}%\\
\address{School of Sciences, South China University of Technology,
Guangzhou 510640, China}
 \email{jing@math.ncsu.edu}

\author[Zhang]{Honglian Zhang$^{\star}$}
\address{Department of Mathematics,
Shanghai University, Shanghai 200444, PR China}
\email{hlzhangmath@shu.edu.cn}
\thanks{$^\star$H.Z., Corresponding Author}

%    General info
\subjclass[2000]{Primary 17B20}
%\date{Version on Aug. 2, 2007}

\keywords{Two-parameter quantum affine algebra,  Young diagram, Fock
space, fermionic realization. }
%%%%%%%%%%%%%%%%%%%%%%%%%%%%%%%%%%%%%%%%%%%%%%%%%%%%%%%%%%%%%%%%%%%%%%%%
%\footnote{Corresponding author.}%
%%%%%%%%%%%%%%%%%%%%%%%%%%%%%%%%%%%%%%%%%%%%%%%%%%%%%%%%%%%%%%%%%%%%%%%%
\begin{abstract}
We construct all fundamental modules for the two parameter quantum
affine algebra of type $A$ using a combinatorial model of Young
diagrams.  In particular we also give a fermionic realization of the
two-parameter quantum affine algebra.

\end{abstract}

\maketitle

\section{ Introduction}

Quantum enveloping algebras are one of the two main examples of
quantum groups introduced by Drinfeld \cite{Dr1} and Jimbo \cite{Jb}
in their study of the Yang-Baxter equation. In the most general
definition the role of the quantum Yang-Baxter equation was given in
the form of universal $R$-matrix. The first example was defined with
the help of Yang-Baxter equation for a $2\times 2$ R-matrix. In the
physics literature and also in the work of Reshetikhin \cite{R} the
dependence on parameters can be more than one variable (see also
\cite{T, J2}).

The research on quantum enveloping algebras with two parameters have
been revitalized by Benkart and Witherspoon in their work on
Drinfeld double construction and Schur-Weyl duality \cite{BW1, BW2,
BW3}. For a review of early history the reader to referred to the
introduction in \cite{BW1}. After then generalizations to other
simple Lie algebras were given by Bergeron, Gao and Hu in
\cite{BGH1, BGH2, BH} and their representations are studied
accordingly. All these work show that the two parameter quantum
groups have similar properties like the usual quantum groups but
offer distinct features in two parameter cases. Recently by
generalizing the vertex representations \cite{FJ, J1} of quantum
affine algebras, Hu, Rosso and the second author \cite{HRZ} have
further introduced a two parameter quantum affine algebra for the
affine type $A$ and also obtained its Drinfeld realization.

It is well-known that quantum affine algebras also admit a fermionic
realization \cite{H}. The fermionic realization has played a
fundamental role in quantum integrable systems and figured
prominently in Kyoto school's work on quantum affine algebras and
their applications to statistical mechanics \cite{DJKMO, JM}. The
crystal basis \cite{Ka1, Ka2} for the quantum general affine algebra
was first constructed with help of Hayashi's fermionic
representation \cite{H}. Our first motivation is to generalize this
construction to construct all level one fundamental representations
for two parameter quantum affine algebra $U_{r,s}(\widehat{\frak
{sl}_n})$. In particular, the fermionic realization is also
obtained, which is two parameter generalization of Misra-Miwa's
realization \cite{MM} and a generalization of Leclerc and Thibon's
version of the combinatorical representation \cite{LT}.

Our second motivation is somewhat more fundamental to justify the
study of two-parametric quantum groups. In the early days the true
meaning of various parameters in quantum groups puzzled some
researchers to question whether the introduction of other parameters
is really necessary. In this paper we will explain the meaning of
two parameters and obtain their combinatorial interpretation, and
show how nicely many pieces of two-parameter quantum groups are
patched together in our fermionic model. Roughly speaking, the two
parameters correspond naturally to the left and right
(multiplication) in our combinatorial model, and one further sees
that the Serre relations are consequence of some combinatorial
properties of our model.

\section{Quantum Affine Algebra $U_{r,s}(\widehat{\frak {sl}_n})$}
\medskip
In this section, we will recall the structure of two-parameter
quantum affine algebra $U_{r,s}(\widehat{\frak {sl}_n})$ defined in
\cite{HRZ}.  Let $\mathbb{K}=\mathbb{Q}(r,s)$ denote the  field of
rational functions in two variables $r$, $s$ ($r\ne \pm s$). For
$n\geq 2$ let $\varepsilon_1, \varepsilon_2, \cdots,
\varepsilon_{n-1}$ be an orthonormal basis of $E=\mathbb{R}^n$ under
the inner product $( \, , )$. Set $I=\{1,\cdots,n-1\}$,
$I_0=\{0\}\cup I$. Then $\Phi=\{\varepsilon_i-\varepsilon_j\mid i
\neq j\in I\}$ is the set of roots for the simple Lie algebra
$\mathfrak{sl}(n)$. We can take
$\Pi=\{\alpha_i=\varepsilon_i-\varepsilon_{i+1}\mid i\in I \}$ as
the basis of simple roots. Let $\delta$ denote the primitive
imaginary root of $\widehat{\frak {sl}_n}$. Take
$\alpha_0=\delta-(\varepsilon_1-\varepsilon_n)$, then
$\Pi'=\{\alpha_i\mid i\in I_0\}$ is a base of simple roots of affine
Lie algebra ${\widehat{\frak{sl}_n}}$.

The following definition is an affinization of the two-parameter
quantum groups for type  $\frak{sl}_n$ (see \cite{BW1}).

\begin{defi}
Let $U=U_{r,s}(\widehat{\frak {sl}_n})$ $(n\geq 2)$ be the unital
associative algebra over $\mathbb{K}$ generated by the elements
$e_j,\, f_j,\, \omega_j^{\pm 1},\, \omega_j'^{\,\pm 1}\, (j\in
I_0)$, $\gamma^{\pm\frac{1}2}$, $\gamma'^{\pm\frac{1}2}$, $
D^{\pm1}$, $D'^{\,\pm1}$, satisfying the following relations:

\medskip
\noindent $(\textrm{A1})$ \
$\gamma^{\pm\frac{1}2},\,\gamma'^{\pm\frac{1}2}$ are central with
$\gamma=\om_\delta$, $\gamma'=\om'_\delta$,  such that
$\omega_i\,\omega_i^{-1}=\omega_i'\,\omega_i'^{\,-1}=1
=DD^{-1}=D'D'^{-1}$, and
\begin{equation*}
\begin{split}[\,\omega_i^{\pm 1},\omega_j^{\,\pm 1}\,]&=[\,\om_i^{\pm1},
D^{\pm1}\,]=[\,\om_j'^{\,\pm1}, D^{\pm1}\,] =[\,\om_i^{\pm1},
D'^{\pm1}\,]=0\\
&=[\,\omega_i^{\pm 1},\omega_j'^{\,\pm 1}\,]=[\,\om_j'^{\,\pm1},
D'^{\pm1}\,]=[D'^{\,\pm1}, D^{\pm1}]=[\,\omega_i'^{\pm
1},\omega_j'^{\,\pm 1}\,].
\end{split}
\end{equation*}
 $(\textrm{A2})$ \ For $\,i,\, j\in I_0$,
\begin{equation*}
\begin{array}{lll}
& D\,e_i\,D^{-1}=r^{\delta_{0i}}\,e_i,\qquad\qquad\qquad\qquad\;
&D\,f_i\,D^{-1}=r^{-\delta_{0i}}\,f_i,\\
&\omega_j\,e_i\,\omega_j^{\,-1}=\langle \omega'_i,\,\omega_j\rangle
\,e_i,\qquad\quad &\omega_j\,f_i\,\omega_j^{\,-1}=\langle
\omega'_i,\,\omega_j\rangle^{-1}\,f_i.
\end{array}
\end{equation*}\\
$(\textrm{A3})$ \  For $\,i,\, j\in I_0$,
\begin{equation*}
\begin{array}{lll}
& D'\,e_i\,D'^{-1}=s^{\delta_{0i}}\,e_i,\qquad\qquad\qquad\quad\ \
&D'\,f_i\,D'^{-1}=s^{-\delta_{0i}}\,f_i,\\
&\omega'_j\,e_i\,\omega'^{\,-1}_j=\langle
\omega'_j,\,\omega_i\rangle^{-1} \,e_i, \qquad\ \
&\omega'_j\,f_i\,\omega'^{\,-1}_j=\langle
\omega'_j,\,\omega_i\rangle\,f_i.
\end{array}
\end{equation*}\\
$(\textrm{A4})$ \ For $\,i,\, j\in I_0$, we have
 $$[\,e_i, f_j\,]=\frac{\delta_{ij}}{r-s}(\omega_i-\omega'_i).$$
$(\textrm{A5})$ \ For $\,i,\,j\in I_0$, but $(i,j)\notin
\{\,(0,n-1), \ (n-1,0)\,\}$ with $a_{ij}=0$, we have
 $$[\,e_i, e_j\,]=0=[\,f_i, f_j\,].$$
$(\textrm{A6})$ \ For $\,i\in I_0$, we have the $(r,s)$-Serre
relations:
\begin{gather*}
e_i^2e_{i+1}-(r+s)\,e_ie_{i+1}e_i+(rs)\,e_{i+1}e_i^2=0,\\
e_ie_{i+1}^2-(r+s)\,e_{i+1}e_ie_{i+1}+(rs)\,e_{i+1}^2e_i=0,\\
e_{n-1}^2e_0-(r+s)\,e_{n-1}e_0e_{n-1}+(rs)\,e_0e_{n-1}^2=0,\\
e_{n-1}e_0^2-(r+s)\,e_0e_{n-1}e_0+(rs)\,e_0^2e_{n-1}=0.
\end{gather*}
$(\textrm{A7})$ \ For $\,i\in I_0$, we have the $(r,s)$-Serre
relations:
\begin{gather*}
f_i^2f_{i+1}-(r^{-1}+s^{-1})\,f_if_{i+1}f_i+(r^{-1}s^{-1})\,f_{i+1}f_i^2=0,\\
f_if_{i+1}^2-(r^{-1}+s^{-1})\,f_{i+1}f_if_{i+1}+(r^{-1}s^{-1})\,f_{i+1}^2f_i=0,\\
f_{n-1}^2f_0-(r^{-1}+s^{-1})\,f_{n-1}f_0f_{n-1}+(r^{-1}s^{-1})\,f_0f_{n-1}^2=0,\\
f_{n-1}f_0^2-(r^{-1}+s^{-1})\,f_0f_{n-1}f_0+(r^{-1}s^{-1})\,f_0^2f_{n-1}=0,
\end{gather*}
where $\langle \omega'_i,\,\omega_j\rangle$ is a skew-dual pairing
defined as follows (more detail see \cite{HRZ}):
$$\langle \omega'_i,\,\omega_j\rangle
=\begin{cases}r^{(\varepsilon_j,\,
\alpha_i)}\,s^{(\varepsilon_{j+1},\, \alpha_i)}, \quad\ \; (i\in
I_0,\, j\in I)\cr r^{-(\varepsilon_{i+1},\,
\alpha_0)}\,s^{(\varepsilon_1,\, \alpha_i)}, \quad (i\in
I_0,\,j=0)\end{cases}$$ From now on, let us write briefly $\langle
\omega'_i,\,\omega_j\rangle=\langle i,\,j\rangle$.
\end{defi}

It can be proved (see \cite{HRZ}) that $U_{r,s}(\widehat{\frak
{sl}_n})$ is a Hopf algebra with the coproduct $\Delta$, the counit
$\vep$ and the antipode $S$ defined below: for $i\in I_0$, we have
\begin{gather*}
\Delta(\gamma^{\pm\frac{1}2})=\gamma^{\pm\frac{1}2}\otimes
\gamma^{\pm\frac{1}2}, \qquad
\Delta((\gamma')^{\,\pm\frac{1}2})=(\gamma')^{\,\pm\frac{1}2}\otimes
(\gamma')^{\,\pm\frac{1}2}, \\
\Delta(D^{\pm1})=D^{\pm1}\otimes D^{\pm1},\qquad
\Delta(D'^{\,\pm1})=D'^{\,\pm1}\otimes D'^{\,\pm1},\\
\Delta(w_i)=w_i\ot w_i, \qquad \Delta(w_i')=w_i'\ot w_i',\\
\Delta(e_i)=e_i\ot 1+w_i\ot e_i, \qquad \Delta(f_i)=f_i\ot w_i'+1\ot
f_i,\\
\vep(e_i)=\vep(f_i)=0,\quad \vep(\gamma^{\pm\frac{1}2})
=\vep((\gamma')^{\,\pm\frac{1}2})=\vep(D^{\pm1})=\vep(D'^{\,\pm1})=\vep(w_i)=\vep(w_i')=1,
\\
S(\gamma^{\pm\frac{1}2})=\gamma^{\mp\frac{1}2},\qquad
S((\gamma')^{\pm\frac{1}2})=(\gamma')^{\mp\frac{1}2},\qquad
S(D^{\pm1})=D^{\mp1},\qquad S(D'^{\,\pm1})=D'^{\,\mp1},\\
S(e_i)=-w_i^{-1}e_i,\qquad S(f_i)=-f_i\,w_i'^{-1},\qquad
S(w_i)=w_i^{-1}, \qquad S(w_i')=w_i'^{-1}.
\end{gather*}

\begin{remark} The algebra $U_{r,s}(\widehat{\frak {sl}_2})$
is isomorphic to the quantum affine algebra $U_{q,
q^{-1}}(\widehat{\frak {sl}_2})$ if set $rs^{-1}=q^2$, see
\cite{HRZ}.
\end{remark}

\begin{remark} We remark that the two sets of generators $\omega_i,
\omega_i'$ follow the original idea of \cite{BW1, BW2} to naturally
blend the second parameter into the relations. Roughly speaking when
one identifies $\omega_i'$ to $\omega_i^{-1}$ the algebra
specializes to the usual quantum affine algebra.
\end{remark}

\section{Fock space representations of $U_{r,s}(\widehat{\frak {sl}_n})$}

In this section we construct a Fock space representation for the
quantum affine algebra $U_{r,s}(\widehat{\frak {sl}_n})$ based on
the fermionic representation of the usual quantum affine algebra.

The Fock space is modeled on the space of partitions. A partition is
a decomposition of a natural number written in nondecreasing order.
Let $\mathcal P(n)$ be the set of partitions of $n$. The generating
function of partitions is given by
$$
\sum_{n=0}^{\infty}|\mathcal
P(n)|q^n=\prod_{n=1}^{\infty}(1-q^n)^{-1}.
$$

For each partition $\lambda=(\lambda_1, \lambda_2, \cdots,
\lambda_l)$ we associate the Young diagram consisting of $n$ nodes
(or boxes) which are stacked in $l$ rows in the 4th quarter of the
xy-coordinate system and aligned at the origin in such a way that
the ith row occupies $\lambda_i$ nodes. The diagonal of the Young
diagram is given by the nodes along the line $y=-x$. Another way to
identify the Young diagram is the following: we specify the path
starting at $(0, -\lambda_1')$ and we move eastward by $\lambda_l$
steps and then we go north by $\lambda_1'-\lambda_2'$ steps and so
on. We say that a node or a box in $\lambda$ sitting at position
$(a, -b)$ if its upper left corner is situated at the point $(a,
-b)$. We will define the residue of any node at $(a, -b)$ to be
$a-b\,( mod \, n)$. For instance the Young diagram
$\lambda=(6,4,4,2,2)$ with residues is given in Figure
\ref{F:young}. For convenience we allow the residues to take values
in $\mathbb Z_n$. So in Figure 1 the residues shown in the lower
left part will be $5, 4, 3$ respectively for the quantum algebra
$U_{r,s}(\widehat{\frak {sl}_6})$.

\begin{figure}[h] \label{F:young}
\setlength{\unitlength}{0.75cm}
\begin{picture}(7,5)
\multiput(0,3)(1,0){7}%
{\line(0,1){1}}
\multiput(0,3)(0,1){2}%
{\line(1,0){6}}      % The above is a row of 5 boxes
\multiput(0,2)(1,0){5}%
{\line(0,1){1}}
\multiput(0,2)(0,1){1}%
{\line(1,0){4}}      % The above is a row of 3 boxes
\multiput(0,1)(1,0){3}%
{\line(0,1){1}}
\multiput(0,1)(0,1){1}%
{\line(1,0){2}}      %
\multiput(0,0)(1,0){3}%
{\line(0,1){1}}
\multiput(0,0)(0,1){1}%
{\line(1,0){2}}      % The above is a row of 1 boxes
\multiput(0.4,3.3)(1, -1){2}{\mbox{$0$}} \multiput(1.4,3.3)(1,
-1){2}{\mbox{$1$}} \multiput(2.4,3.3)(1,-1){2}{\mbox{$2$}}
\multiput(3.4,3.3)(1,0){1}{\mbox{$3$}}
\multiput(4.4,3.3)(1,0){1}{\mbox{$4$}}\put(5.4,3.3){\mbox{$5$}}
\multiput(0.1, 2.3)(1,-1){2}{\mbox{${-1}$}} \multiput(0.1,
1.3)(1,-1){2}{\mbox{$-2$}} \put(0.1, 0.3){\mbox{$-3$}}
\end{picture}\par
\caption{Young diagram $(64^22^2)$}
\end{figure}

A node $\gamma$ or box (convex corner) is called {\it removable} if
$\lambda-\gamma$ is still a Young diagram. A contract corner is
called {\it intent} if one node or box can be added at the corner to
form another Young diagram. By abusing the terminology we will call
a node $(a, b)$ of the border rim of the diagram of $\lambda$ {\it
indent} if the box $(a-1, b-1)$, $(a+1, b)$ or $(a, b-1)$ can be
added to $\lambda$ to form a new Young diagram. Note that the new
diagram's border does not contain the node except the node is a
starting node or a final node of the border.

For example in Figure 1 the nodes at $(1, -1)$, $(3, 0), (5, 0), (0,
-3)$ are indent nodes except that the nodes at $(5, 0)$, $(3, -1)$,
$(1, -3)$ are removable nodes. By convention the trivial Young
diagram $\phi$ has an imaginary indent node at $(-1, 1)$. For this
reason we sometimes say that a point at $(a, -b)$ is removable or
indented. In this sense the origin for the trivial Young diagram is
indented.

For each $i\in \mathbb Z_n$ we define the Young diagram $|\lambda,
i\rangle$ as the diagram that assigns  the residues $i+a-b \, ( mod
\, n)$ to each node $(a, -b)$ in $\lambda$. The previous example is
the Young diagram for $i=0$. The following Figure 2 shows the
residue of $(a, -b)$ for the  configuration of $|\lambda, i\rangle$.

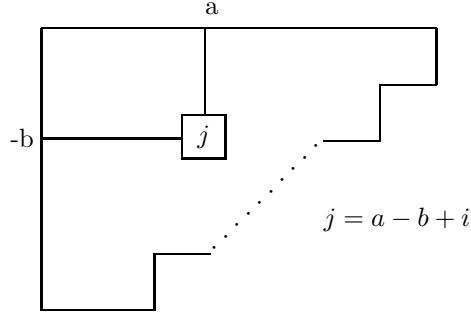
\begin{figure}[h] \label{F:young2}
\setlength{\unitlength}{0.75cm}
\begin{picture}(7,7)
\put(0, 5){\line(1,0){7}} \put(0,0){\line(0,1){5}}
\put(7,4){\line(0,1){1}}\put(7, 4){\line(-1,0){1}}
\put(6,4){\line(0,-1){1}}\put(6, 3){\line(-1,0){1}}

\put(0,0){\line(1,0){2}}\put(2,0){\line(0,1){1}}
\put(2,1){\line(1,0){1}} \multiput(3,1)(0.2,0.2){10}{\mbox{$\cdot$}}

\put(2.5,2.7){\framebox(0.75,0.75){$ j$}} \put(0,
3.05){\line(1,0){2.5}}\put(2.9, 5){\line(0,-1){1.55}} \put(-0.55,
2.9){\mbox{-b}}\put(2.9, 5.25){\mbox{a}}
\put(5,1.5){\mbox{$j=a-b+i$}}

\end{picture}\par
\caption{Residue for Young diagram $|\lambda, i\rangle$}
\end{figure}

Let $K=\mathbb Q(r, s)$ be the field of rational functions in $r$
and $s$. For each $i\in \mathbb Z_n$ we define the Fock space
$\mathcal F_i$ to be
\begin{equation}
\mathcal F_i=\bigoplus_{\lambda\in \mathcal P} K|\lambda, i\rangle.
\end{equation}

Let $\lambda$ and $\mu$ be two partitions such that $\mu$ is
obtained from $\lambda$ by adding a node $\gamma$ of residue $i$, or
$\mu/\lambda$ consists of a single $i$-node $\gamma$. Let
$I_i(\lambda)$ be the set of indent $i$-nodes in the boundary of
$\lambda$,
 and $R_i(\lambda)$ be the set of its removable $i$-nodes in the boundary of
 $\lambda$.
We further let $I_i^{(l)}(\lambda,\mu)$ (resp.
$R_i^{(l)}(\lambda,\mu)$) be the set of indent $i$-nodes (resp. of
removable $i$-nodes) situated to the left of $\gamma$ ($\gamma$ not
included) in the boundary of $\mu$, and similarly, let
$I_i^{(r)}(\lambda,\mu)$ and $R_i^{(r)}(\lambda,\mu)$ be  the
corresponding set of $i$-nodes located not to the left of $\gamma$
(i.e. to the right of $\gamma$ and $\gamma$ itself) on the boundary
of $\lambda$. We also denote by $I^{(0)}(\lambda)$ ( resp.
$R^{(0)}(\lambda)$ ) the total number of indent (resp. removable)
$0$-nodes in diagram $\lambda$. We denote by
$|I_i^{(r)}(\lambda,\mu)|$ and $|R_i^{(r)}(\lambda,\mu)|$ the number
of the set $I_i^{(r)}(\lambda,\mu)$ and
$R_i^{(r)}(\lambda,\mu)$ respectively.\\[6pt]

For a finite set $S$, we use $|S|$ to denote its cardinality. We
define the action of the simple generators as follows.
\begin{eqnarray*}
 && f_i~|~\lambda~ \rangle = \sum_{\mu}
r^{|I_i^{(r)}(\lambda,~\mu)|}~s^{|R_i^{(r)}(\lambda,~\mu)|}~| ~\mu
~\rangle  ,\\
 && e_i~|~\lambda~ \rangle = \sum_{\mu}
r^{|R_i^{(l)}(\mu,~ \lambda)|}~s^{|I_i^{(l)}(\mu, ~\lambda)|}~| ~\mu
~\rangle ,\\
 && \omega_i~|~\lambda~ \rangle = r^{|I_i(\lambda)|}~s^{|R_i(\lambda)|}~|
 ~\lambda ~\rangle ,\\[6pt]
 &&\omega'_i~|~\lambda~ \rangle = r^{|R_i(\lambda)|}~s^{|I_i(\lambda)|}~|
 ~\lambda ~\rangle,\\[6pt]
&& D~|~\lambda~ \rangle = r^{- |I^0(\lambda)|}~s^{-|R^0(\lambda)|}~|
 ~\lambda ~\rangle,\\[6pt]
 && D'_i~|~\lambda~ \rangle = r^{-|R^0(\lambda)|}~s^{-|I^0(\lambda)|}~|
 ~\lambda ~\rangle,
\end{eqnarray*}
where the first sum runs through all $\mu$ such that $\mu/\lambda$
is a $i$-node, and the second sum runs through all $\mu$ such that
$\lambda/\mu$ is a $i$-node. We remark that all the sums are
obviously finite since it runs through the cells of the Young
diagram. We also note that when $r=s^{-1}=q$, the above action
reduces to that of the quantum affine algebra $U_q(\widehat{sl}_n)$
formulated in \cite{LT}.

\begin{theorem} For each $i\in\mathbb Z_n$ the Fock space
$\mathcal F_i$ is a level one representation of $\mathcal
{U}_{r,s}(\hat{sl}_n) $ and contains $V(\Lambda_i)$ as a submodule
and the highest weight vector is $|\phi, i>$, where $\phi$ is the
empty diagram.
\end{theorem}

\noindent{\bf Proof.} First of all we note the statement about the
submodule is obtained by computing the action of the Heisenberg
subalgebra. By symmetry it is enough to show the statement for the
basic module $V(\Lambda_0)$, so we will drop the second index $i=0$
in $|\lambda, 0\rangle$. Namely, we take $i=0$ and the residue for
the node $(a, -b)$ is $a-b$.

(i) From the construction we first have
\begin{eqnarray*}
  \omega_j e_i~|~\lambda~\rangle &= &\omega_j \sum_{\gamma \in R_i(\lambda)}r^{|R_i^{(l)}(\lambda-\gamma,~\lambda)|}
   ~s^{|I_i^{(l)}(\lambda-\gamma,~\lambda)|} ~|~\lambda-\gamma~\rangle\\
   &=& \sum_{\gamma \in R_i(\lambda)}r^{|R_i^{(l)}(\lambda-\gamma,~\lambda)|}
   ~s^{|I_i^{(l)}(\lambda-\gamma,~\lambda)|}r^{|I_j(\lambda-\gamma)|}~s^{|R_j(\lambda-\gamma)|}~
   |~\lambda-\gamma~\rangle,
\end{eqnarray*}
and similarly
\begin{equation*}
   e_i \omega_j~|~\lambda~\rangle=\sum_{\gamma \in R_i(\lambda)}r^{|R_i^{(l)}(\lambda-\gamma,~\lambda)|}
   ~s^{|I_i^{(l)}(\lambda-\gamma,~\lambda)|} ~ r^{|I_j(\lambda)|}
   ~s^{|R_j(\lambda)|}~|~\lambda-\gamma~\rangle.
\end{equation*}
Since for $\gamma\in R_i(\lambda)$,
\begin{eqnarray} \label{1}
 |I_i(\lambda-\gamma)|= |I_i(\lambda)|+1,~~~|R_i(\lambda-\gamma)|
=|R_i(\lambda)|-1,
\end{eqnarray}
It follows immediately from (\ref{1})  that
$$\omega_i e_i = rs^{-1} e_i \omega_i = \langle i,\,i\rangle e_i \omega_i.$$

On the other hand it follows from $\gamma \in R_i(\lambda)~~~(i\neq
j) $ that,
\begin{eqnarray}
 &&|I_j(\lambda-\gamma)|= |I_j(\lambda)|+ \begin{cases}(\varepsilon_j,\,
\alpha_i), \quad\ \; (i\in I_0,\, j\in I)\vspace{6pt} \cr
{-(\varepsilon_{i+1},\, \alpha_0)}, \quad (i\in
I_0,\,j=0)\end{cases}\\[6pt]
 &&|R_j(\lambda-\gamma)|
=|R_j(\lambda)|+\begin{cases}{(\varepsilon_{j+1},\, \alpha_i)},
\quad\ \; (i\in I_0,\, j\in I) \vspace{6pt} \cr {(\varepsilon_1,\,
\alpha_i)}, \quad (i\in I_0,\,j=0).\end{cases}
\end{eqnarray}
These combinatorial identities imply that
$$\omega_j e_i =  <i,~j> e_i \omega_j.$$

(ii) To check the commutation relation (A4), we consider the actions
of $e_i$ and $f_i$:
\begin{eqnarray*}
 && e_i f_i ~|~\lambda~\rangle =  e_i \sum_{\gamma \in I_i(\lambda)}r^{|I_i^{(r)}(\lambda,~\lambda+\gamma)|}
   ~s^{|R_i^{(r)}(\lambda,~\lambda+\gamma)|} ~|~\lambda+\gamma~\rangle\\
  &&  = \sum_{ {\gamma \in I_i(\lambda),}\atop {\gamma' \in R_i(\lambda+\gamma)}}\begin{array}{l}r^{|R_i^{(r)}(\lambda,~\lambda+\gamma)|}
   ~s^{|R_i^{(r)}(\lambda,~\lambda+\gamma)|}\\ \times~  r^{|R_i^{(l)}(\lambda+\gamma-\gamma',\lambda+\gamma)|}
   ~s^{|I_i^{(l)}(\lambda+\gamma-\gamma',\lambda+\gamma)|}
   |\lambda+\gamma-\gamma'\rangle \end{array}\\
&&  = \sum_{ {\gamma \in I_i(\lambda),}\atop {\gamma' \in
R_i(\lambda)\cup I_i(\lambda)}} \begin{array}{l}
r^{|I_i^{(r)}(\lambda,~\lambda+\gamma)|+|R_i^{(l)}(\lambda+\gamma-\gamma',\lambda+\gamma)|}
  \\ ~\times s^{|R_i^{(r)}(\lambda,~\lambda+\gamma)|+|I_i^{(l)}(\lambda+\gamma-\gamma',\lambda+\gamma)|}~
   ~|~\lambda+\gamma-\gamma'\rangle \end{array}\\
&&  = \sum_{\gamma \in I_i(\lambda)}r^{|I_i^{(r)}(\lambda,
\lambda+\gamma)|+|R_i^{(l)}(\lambda,\lambda+\gamma)|~
s^{|R_i^{(r)}}(\lambda,\lambda+\gamma)|+|I_i^{(l)}(\lambda,\lambda+\gamma)|} ~~|~\lambda\rangle \\
&&~~~~~ +\sum_{ {\gamma \in I_i(\lambda),}\atop {\gamma' \in
R_i(\lambda)}}
\begin{array}{l}r^{|I_i^{(r)}(\lambda,~\lambda+\gamma)|+|R_i^{(l)}(\lambda+\gamma-\gamma',\lambda+\gamma)|}
   \\~ \times s^{|R_i^{(r)}(\lambda,~\lambda+\gamma)|+|I_i^{(l)}(\lambda+\gamma-\gamma',\lambda+\gamma)|}~
   |\lambda+\gamma-\gamma'\rangle\end{array},
\end{eqnarray*}
where we used the result: $R_i(\lambda+\gamma)=R_i(\lambda)+\{
\gamma \}$.

Reversing the order of the product we have,
\begin{eqnarray*}
 && f_i e_i ~|~\lambda~\rangle
 %=  f_i \sum_{\gamma' \in R_i(\lambda)}r^{|R_i^{(l)}(\lambda-\gamma',~\lambda)|}
 %  ~s^{|I_i^{(l)}(\lambda-\gamma',~\lambda)|} ~|~\lambda-\gamma'~\rangle\\&&
  = \sum_{ {\gamma' \in R_i(\lambda),}\atop {\gamma \in I_i(\lambda-\gamma')}} \begin{array}{l}
  r^{|R_i^{(l)}(\lambda-\gamma',
  ~\lambda)|}
   ~s^{|I_i^{(r)}(\lambda-\gamma',~\lambda)|}\\  ~\times~ r^{|I_i^{(r)}(\lambda-\gamma',\lambda-\gamma'+\gamma)|}
   ~s^{|R_i^{(r)}(\lambda-\gamma',\lambda-\gamma'+\gamma)|}~
   |~\lambda-\gamma'+\gamma~\rangle \end{array}\\
&&  = \sum_{\gamma' \in R_i(\lambda)}r^{|R_i^{(l)}(\lambda-\gamma',
\lambda)|+|I_i^{(r)}(\lambda-\gamma',\lambda)|~
s^{|I_i^{(l)}}(\lambda-\gamma',\lambda)|+|R_i^{(r)}(\lambda-\gamma',\lambda)|} ~|~\lambda~\rangle \\
&&~~~~~~~ +\sum_{ {\gamma' \in R_i(\lambda),} \atop {\gamma' \in
I_i(\lambda)}} \begin{array}{l}
r^{|R_i^{(l)}(\lambda-\gamma',~\lambda)|+|I_i^{(r)}(\lambda-\gamma',\lambda-\gamma'+\gamma)|}
   \\ ~\times~s^{|I_i^{(r\l)}(\lambda-\gamma',~\lambda)|+|R_i^{(r)}(\lambda-\gamma',\lambda-\gamma'+\gamma)|}~
   |~\lambda-\gamma'+\gamma~\rangle\end{array},\\
\end{eqnarray*}
The following fact is easily verified.

\noindent{\bf Claim A}\, For all $\gamma \in I_i(\lambda),~~ \gamma'
\in
R_i(\lambda)$ \\
\begin{eqnarray}
 {\label{2}} && I_i^{(r)}(\lambda,\lambda+\gamma) -I_i^{(r)}(\lambda-\gamma',
 \lambda-\gamma'+\gamma)   \\
&&\hskip2cm = R_i^{(l)}(\lambda-\gamma')- R_i^{(l)}(\lambda
 +\gamma-\gamma',\lambda+\gamma);\nonumber\\
{\label{3}} &&
I_i^{(l)}(\lambda-\gamma',\lambda)-I_i^{(l)}(\lambda+\gamma-\gamma',\lambda+
\gamma)\\
&&\hskip2cm =
 R_i^{(r)}(\lambda,\lambda+\gamma)-R_i^{(r)}(\lambda-\gamma',\lambda-\gamma'+\gamma).\nonumber
\end{eqnarray}
Combining the above two expressions in Claim A, we get that
\begin{eqnarray*}
& & [\,e_i,\, f_i\,] ~|~\lambda~\rangle \\
& =&
  \sum_{\gamma \in I_i(\lambda)}r^{|I_i^{(r)}(\lambda,\,~\lambda+\gamma)|+|R_i^{(l)}(\lambda,\,~\lambda+\gamma)|}
   ~s^{|I_i^{(l)}(\lambda,\,~\lambda+\gamma)|+|R_i^{(r)}(\lambda,\,~\lambda+\gamma)|}|~\lambda~\rangle\\
  && - \sum_{\gamma' \in
R_i(\lambda)}r^{|I_i^{(r)}(\lambda-\gamma',\,~\lambda)|+|R_i^{(l)}(\lambda-\gamma',\,~\lambda)|}
   ~s^{|I_i^{(l)}(\lambda-\gamma',\,~\lambda)|+|R_i^{(r)}(\lambda-\gamma',\,~\lambda)|}|~\lambda~\rangle
  \end{eqnarray*}
The following Claim B is important for the further deduction.\\

\noindent{\bf Claim B} \,\, For all $\gamma \in I_i(\lambda),~~
\gamma' \in
R_i(\lambda)$ \\
\begin{eqnarray}
 && |I_i(\lambda-\gamma')|=|I_i^{(r)}(\lambda-\gamma',\,\lambda)| +|I_i^{(l)}(\lambda-\gamma',
 \lambda)|+1;\\[3pt]
 && |R_i(\lambda+\gamma)|=|R_i^{(r)}(\lambda,\,\lambda+\gamma)| +|R_i^{(l)}(\lambda,\,\lambda+\gamma)|+1;\\[3pt]
&& |I_i(\lambda)|=|I_i^{(r)}(\lambda,\,\lambda+\gamma)| +|I_i^{(l)}(\lambda,\,\lambda+\gamma)|+1;\\[3pt]
 && |R_i(\lambda)|=|R_i^{(r)}(\lambda-\gamma',\,\lambda)| +|R_i^{(l)}(\lambda-\gamma',\,\lambda)|+1;
\end{eqnarray}

\medskip

We note that the coefficient is given by
\begin{eqnarray*}
  && \sum_{\gamma \in I_i(\lambda)}r^{|I_i^{(r)}(\lambda,\,~\lambda+\gamma)|+|R_i^{(l)}(\lambda,\,~\lambda+\gamma)|}
   ~s^{|I_i^{(l)}(\lambda,\,~\lambda+\gamma)|+|R_i^{(r)}(\lambda,\,~\lambda+\gamma)|}\\[3pt]
  &-&  \sum_{\gamma' \in R_i(\lambda)}r^{|I_i^{(r)}(\lambda-\gamma',\,~\lambda)|+|R_i^{(l)}(\lambda-\gamma',\,~\lambda)|}
   ~s^{|I_i^{(l)}(\lambda-\gamma',\,~\lambda)|+|R_i^{(r)}(\lambda-\gamma',\,~\lambda)|}\\[3pt]
&=&\sum_{\gamma \in
I_i(\lambda)}r^{|I_i^{(r)}(\lambda,\,~\lambda+\gamma)|-|R_i^{(r)}(\lambda,\,~\lambda+\gamma)|+|R_i(\lambda+\gamma)|-1}
   ~s^{|I_i(\lambda)|-|I_i^{(r)}(\lambda,\,~\lambda+\gamma)|+|R_i^{(r)}(\lambda,\,~\lambda+\gamma)|-1}\\
&-& \sum_{\gamma' \in
R_i(\lambda)}r^{|I_i^{(r)}(\lambda-\gamma',\,~\lambda)|-|R_i^{(r)}(\lambda-\gamma',\,~\lambda)|+|R_i(\lambda)|-1}
   ~s^{|I_i(\lambda-\gamma')|-|I_i^{(r)}(\lambda-\gamma',\,~\lambda)|+|R_i^{(r)}(\lambda-\gamma',\,~\lambda)|-1}
\end{eqnarray*}
where we used  the relations (3.8) and (3.9) in the first term and
the relations (3.7) and (3.10) in the second term. Then we have
\begin{eqnarray*}
\text{The action of}\quad [e_i, \,f_i]&=&\sum_{\gamma \in
I_i(\lambda)}r^{|R_i(\lambda+\gamma)|}s^{|I_i(\lambda)|-2}
(rs^{-1})^{|I_i^{(r)}(\lambda,\,~\lambda+\gamma)|-|R_i^{(r)}(\lambda,\,~\lambda+\gamma)|-1}\\[3pt]
&-& \sum_{\gamma' \in
R_i(\lambda)}r^{|R_i(\lambda)|+1}s^{|I_i(\lambda-\gamma')|-3}
(rs^{-1})^{|I_i^{(r)}(\lambda-\gamma',\,\lambda)|-|R_i^{(r)}(\lambda-\gamma',\,\lambda)|-2}\\[3pt]
&=&r^{|R_i(\lambda)|+1}s^{|I_i(\lambda)|-2}\Big(\sum_{\gamma \in
I_i(\lambda)}
(rs^{-1})^{|I_i^{(r)}(\lambda,\,~\lambda+\gamma)|-|R_i^{(r)}(\lambda,\,~\lambda+\gamma)|-1}\\[3pt]
&&\hskip3cm - \sum_{\gamma' \in R_i(\lambda)}
(rs^{-1})^{|I_i^{(r)}(\lambda-\gamma',\,\lambda)|-|R_i^{(r)}(\lambda-\gamma',\,\lambda)|-2}\Big)\\[3pt]
&=&r^{|R_i(\lambda)|+1}s^{|I_i(\lambda)|-2}\Big(
(rs^{-1})^{|I_i(\lambda)|-|R_i(\lambda)|-2}+(rs^{-1})^{|I_i(\lambda)|-|R_i(\lambda)|-3}\\[3pt]
&&\hskip3.2cm +\cdots +
(rs^{-1})^{|I_i(\lambda)|-|R_i(\lambda)|-|I_i(\lambda)|+|R_i(\lambda)|-1}\Big)\\[3pt]
%&=&r^{|I_i(\lambda)|-1}s^{|R_i(\lambda)|}+r^{|I_i(\lambda)|-2}s^{|R_i(\lambda)|+1}
% +\cdots + r^{|R_i(\lambda)|}s^{|I_i(\lambda)|-1}\\[3pt]
&=&\frac{r^{|I_i(\lambda)|}s^{|R_i(\lambda)|}-r^{|R_i(\lambda)|}s^{|I_i(\lambda)|}}{r-s}
=\text{the action of} \quad \frac{\omega_i-\omega'_i}{r-s}
\end{eqnarray*}
%Therefore we have checked the commutative relation:
%$$[e_i, \,f_i]=\frac{\omega_i-\omega'_i}{r-s}.$$

Finally, we will check the $(r, s)-$ Serre relation:

$$f_{i+1}f_i^2-(r+s)\,f_if_{i+1}f_i+rs\,f_i^2f_{i+1}=0.$$

It follows from the definition that
 \begin{eqnarray*}
&& f_{i+1}f_i^2~|~\lambda~\rangle \\
&=& \sum_{\mbox{\tiny $\begin{array}{c}\gamma_1 \in
 I_i(\lambda) \\ \gamma_2\in I_i(\lambda+\gamma_1) \\ \gamma_3\in I_{i+1}(\lambda
 +\gamma_1+\gamma_2)\end{array}$}}
 \begin{array}{l}
 r^{|I_i^{(r)}(\lambda,\, \lambda+\gamma_1)|+|I_i^{(r)}(\lambda+\gamma_1,\, \lambda+\gamma_1+\gamma_2)|
 +|I_{i+1}^{(r)}(\lambda +\gamma_1+\gamma_2,\, \lambda
 +\gamma_1+\gamma_2+\gamma_3)|}\\[3pt]
  s^{|R_i^{(r)}(\lambda,\, \lambda+\gamma_1)|+
  |R_i^{(r)}(\lambda+\gamma_1,\,
 \lambda +\gamma_1+\gamma_2)|+
 |R_{i+1}^{(r)}(\lambda+\gamma_1+\gamma_2,\,
 \lambda+\gamma_1+\gamma_2+\gamma_3)|}\\[3pt]
 ~~\,\,|~~\lambda+\gamma_1+\gamma_2+\gamma_3\rangle \end{array}
\end{eqnarray*}

For simplicity, we write $$I_{i+1}(\lambda
 +\gamma_1+\gamma_2)=I_{i+1}(\lambda)\cup \bigg(I_{i+1}(\lambda
 +\gamma_1)-I_{i+1}(\lambda)\bigg)\cup \bigg(I_{i+1}(\lambda
 +\gamma_2)-I_{i+1}(\lambda)\bigg),$$
which is derived since $I_{i+1}(\lambda+\gamma_1+\gamma_2)=
I_{i+1}(\lambda)\cup I_{i+1}(\lambda+\gamma_1)\cup
I_{i+1}(\lambda+\gamma_2)$.

  Then the coefficient of the above expression becomes
\begin{eqnarray*}
& &   \sum_{\mbox{\tiny $\begin{array}{c}\gamma_1 \in
 I_i(\lambda)\\
 \gamma_2\in I_i(\lambda+\gamma_1) \\ \gamma_3\in I_{i+1}(\lambda)\end{array}$}}
 \begin{array}{l}r^{|I_i^{(r)}(\lambda,\,\lambda+\gamma_1)|+|I_i^{(r)}(\lambda+\gamma_1,\, \lambda+\gamma_1+\gamma_2)|
 +|I_{i+1}^{(r)}(\lambda,\, \lambda+\gamma_3)|}\\[6pt]
  \times s^{|R_i^{(r)}(\lambda,\lambda+\gamma_1)|+
  |R_i^{(r)}(\lambda+\gamma_1,\,
 \lambda +\gamma_1+\gamma_2)|+ |R_{i+1}^{(r)}(\lambda,\,
 \lambda+\gamma_3)|}\end{array}\\
 && {\displaystyle \hskip 5mm  +  \sum_{\mbox{\tiny $\begin{array}{c}\gamma_1 \in
 I_i(\lambda)\\
 \gamma_2\in I_i(\lambda+\gamma_1) \\ \gamma_3\in I_{i+1}(\lambda
 +\gamma_1) - I_{i+1}(\lambda)\end{array}$}}
 \begin{array}{l}
 r^{|I_i^{(r)}(\lambda,\, \lambda+\gamma_1)|+|I_i^{(r)}(\lambda+\gamma_1,\, \lambda+\gamma_1+\gamma_2)|
 +|I_{i+1}^{(r)}(\lambda +\gamma_1,\, \lambda
 +\gamma_1+\gamma_3)|+1}\\[6pt]
  \times s^{|R_i^{(r)}(\lambda,\, \lambda+\gamma_1)|+
  |R_i^{(r)}(\lambda+\gamma_1,\,
 \lambda +\gamma_1+\gamma_2)|+ |R_{i+1}^{(r)}(\lambda+\gamma_1,\,
 \lambda+\gamma_1+\gamma_3)|}\end{array}}\\
&& {\displaystyle \hskip 5mm   +  \sum_{\mbox{\tiny
$\begin{array}{c}\gamma_1 \in
 I_i(\lambda)\\
 \gamma_2\in I_i(\lambda+\gamma_1) \\ \gamma_3\in
 I_{i+1}(\lambda+\gamma_2)- I_{i+1}(\lambda)\end{array}$}}
 \begin{array}{l} r^{|I_i^{(r)}(\lambda,\, \lambda+\gamma_1)|+|I_i^{(r)}(\lambda+\gamma_1,\, \lambda+\gamma_1+\gamma_2)|
 +|I_{i+1}^{(r)}(\lambda+\gamma_2,\, \lambda+\gamma_2+\gamma_3)|+1}\\[6pt]
  \times s^{|R_i^{(r)}(\lambda,\, \lambda+\gamma_1)|+
  |R_i^{(r)}(\lambda+\gamma_1,\,
 \lambda +\gamma_1+\gamma_2)|+ |R_{i+1}^{(r)}(\lambda+\gamma_2,\,
 \lambda+\gamma_2+\gamma_3)|}\end{array}}\\
\end{eqnarray*}

Furthermore we get
\begin{eqnarray*}
 & & f_i f_{i+1}f_i~|~\lambda~\rangle \\
  &&  = \sum_{\mbox{\tiny $\begin{array}{c}\gamma_1 \in
 I_i(\lambda)\\
 \gamma_3\in I_{i+1}(\lambda+\gamma_1) \\ \gamma_2\in I_i(\lambda
 +\gamma_1+\gamma_3)\end{array}$}}
 \begin{array}{l}
 r^{|I_i^{(r)}(\lambda,\, \lambda+\gamma_1)|+|I_{i+1}^{(r)}(\lambda+\gamma_1,\, \lambda+\gamma_1+\gamma_3)|
 +|I_i^{(r)}(\lambda +\gamma_1+\gamma_3,\, \lambda
 +\gamma_1+\gamma_3+\gamma_2)|}\\[3pt]
   s^{|R_i^{(r)}(\lambda,\, \lambda+\gamma_1)|+
   |R_{i+1}^{(r)}(\lambda+\gamma_1,\,
 \lambda +\gamma_1+\gamma_3)|+
 |R_i^{(r)}(\lambda+\gamma_1+\gamma_3,\,
 \lambda+\gamma_1+\gamma_3+\gamma_2)|}\\[3pt]~~\,\,|~~\lambda+\gamma_1+\gamma_2+\gamma_3\rangle\end{array}
\end{eqnarray*}
Similarly, using $$I_{i+1}(\lambda
 +\gamma_1)=I_{i+1}(\lambda)\cup \bigg(I_{i+1}(\lambda
 +\gamma_1)-I_{i+1}(\lambda)\bigg),$$
the coefficient of the second expression becomes
\begin{eqnarray*}
& & \sum_{\mbox{\tiny $\begin{array}{c}\gamma_1 \in
 I_i(\lambda)\\
 \gamma_3\in I_{i+1}(\lambda) \\ \gamma_2\in I_i(\lambda+\gamma_1)\end{array}$}}
 \begin{array}{l}r^{|I_i^{(r)}(\lambda,\, \lambda+\gamma_1)|+|I_{i+1}^{(r)}(\lambda,\, \lambda+\gamma_3)|
 +|I_i^{(r)}(\lambda+\gamma_1,\, \lambda+\gamma_1+\gamma_2)|}\\[6pt]
  s^{|R_i^{(r)}(\lambda,\, \lambda+\gamma_1)|+
  |R_{i+1}^{(r)}(\lambda,\,
 \lambda+\gamma_3)|+ |R_i^{(r)}(\lambda+\gamma_1,\,
 \lambda +\gamma_1+\gamma_2)|-1}\end{array}\\
 && {\displaystyle \hskip 5mm  +  \sum_{\mbox{\tiny $\begin{array}{c}\gamma_1 \in
 I_i(\lambda)\\
 \gamma_3\in I_{i+1}(\lambda+\gamma_1)-I_{i+1}(\lambda) \\ \gamma_2\in I_i(\lambda
 +\gamma_1)\end{array}$}}
 \begin{array}{l}
 r^{|I_i^{(r)}(\lambda,\, \lambda+\gamma_1)|+|I_{i+1}^{(r)}(\lambda+\gamma_1,\, \lambda+\gamma_1+\gamma_3)|
 +|I_i^{(r)}(\lambda +\gamma_1,\, \lambda
 +\gamma_1+\gamma_2)|+1}\\[6pt]
  s^{|R_i^{(r)}(\lambda,\, \lambda+\gamma_1)|+
  |R_{i+1}^{(r)}(\lambda+\gamma_1,\,
 \lambda +\gamma_1+\gamma_3)|+ |R_i^{(r)}(\lambda+\gamma_1,\,
 \lambda+\gamma_1+\gamma_2)|-1}\end{array}}
\end{eqnarray*}
Finally, using the definition again one has
\begin{eqnarray*}
 & & f_i^2 f_{i+1}~|~\lambda~\rangle \\
  &&  = \sum_{\mbox{\tiny $\begin{array}{c}\gamma_3 \in
 I_{i+1}(\lambda)\\
 \gamma_1\in I_i(\lambda+\gamma_3) \\ \gamma_2\in I_i(\lambda
 +\gamma_1+\gamma_3)\end{array}$}}
 \begin{array}{l}
 r^{|I_{i+1}^{(r)}(\lambda,\, \lambda+\gamma_3)|+|I_i^{(r)}(\lambda+\gamma_3,\, \lambda+\gamma_3+\gamma_1)|
 +|I_i^{(r)}(\lambda +\gamma_1+\gamma_3,\, \lambda
 +\gamma_1+\gamma_3+\gamma_2)|}\\[3pt]
 s^{|R_{i+1}^{(r)}(\lambda,\, \lambda+\gamma_1)|+
 |R_i^{(r)}(\lambda+\gamma_3,\,
 \lambda +\gamma_3+\gamma_1)|+
 |R_i^{(r)}(\lambda+\gamma_1+\gamma_3,\,
 \lambda+\gamma_1+\gamma_2+\gamma_3)|}\\[3pt]~~\,\,|~~\lambda+\gamma_1+\gamma_2+\gamma_3\rangle\end{array}
\end{eqnarray*}
Since for $\gamma_3\in I_{i+1}(\lambda)$, one gets
$$I_i(\lambda+\gamma_3)=I_i(\lambda), \qquad I_i(\lambda
 +\gamma_1+\gamma_3)=I_i(\lambda
 +\gamma_1).$$
 Thus the coefficient of the third expression becomes
\begin{eqnarray*}
 && \sum_{\mbox{\tiny $\begin{array}{c}\gamma_3 \in
 I_{i+1}(\lambda)\\
 \gamma_1\in I_i(\lambda) \\ \gamma_2\in I_i(\lambda
 +\gamma_1)\end{array}$}}
 \begin{array}{l}
 r^{|I_{i+1}^{(r)}(\lambda,\, \lambda+\gamma_3)|+|I_i^{(r)}(\lambda,\, \lambda+\gamma_1)|
 +|I_i^{(r)}(\lambda +\gamma_1,\, \lambda
 +\gamma_1+\gamma_2)|}\\[3pt]
  \times s^{|R_{i+1}^{(r)}(\lambda,\, \lambda+\gamma_3)|+
  |R_i^{(r)}(\lambda,\,
 \lambda +\gamma_1)|-1+ |R_i^{(r)}(\lambda+\gamma_1,\,
 \lambda+\gamma_1+\gamma_2)|-1}\end{array}\\
 \end{eqnarray*}

Combining the above three coefficients, we get the required result:
$$f_{i+1}f_i^2-(r+s)\,f_if_{i+1}f_i+rs\,f_i^2f_{i+1}|\lambda \rangle=0,$$

where we have used the following fact.\\[6pt]

\noindent{\bf Claim C}\, \, For $\gamma_2\in I_i(\lambda+\gamma_1)$
and  $\gamma_3\in I_{i+1}(\lambda+\gamma_1+\gamma_2)$,  it follows
that,
 \begin{eqnarray*}
 && | I_{i+1}^{(r)} (\lambda+\gamma_2, ~ \lambda +\gamma_2+\gamma_3)| =
     |I_{i+1}^{(r)} (\lambda+\gamma_1, ~ \lambda +\gamma_1+\gamma_3)|+1,\\[6pt]&&
     |R_{i+1}^{(r)} (\lambda+\gamma_2, ~ \lambda +\gamma_2+\gamma_3)| =
     |R_{i+1}^{(r)} (\lambda+\gamma_1, ~ \lambda
     +\gamma_1+\gamma_3)|-1.
 \end{eqnarray*}
\medskip
This completes the proof of Theorem 3.1.

\vskip30pt \centerline{\bf ACKNOWLEDGMENT}

\bigskip

N. Jing would like to thank the support of NSA grant and NSFC (No.
10728102). H. Zhang would like to thank the support of NSFC (No.
10801094)and Shanghai Leading Academic Discipline Project(No.
J50101).
\bigskip

\bibliographystyle{amsalpha}

\end{document}